\documentclass[a4paper,twocolumn]{article}
\usepackage{fullpage}
\usepackage[latin1]{inputenc}
\usepackage[T1]{fontenc}
\usepackage{palatcm}
\usepackage{amsfonts,amssymb,stmaryrd,amsmath,bbm}
\usepackage{graphicx,epsfig,color,ifpdf}

\newcommand{\ZZ}{\mathbb Z}
\newcommand{\un}{\mathbbm{1}}

\let\epsilon=\varepsilon

\title{\huge\bf Percolation theory}
\author{\bf Vincent  \textsc{Beffara} --- CNRS, UMPA, ENS Lyon\\
\bf Vladas \textsc{Sidoravicius} --- IMPA, Rio de Janeiro}
\date{Printed \today}

\newtheorem{thm}{Theorem}

\begin{document}
\maketitle


\section*{Introduction}

Percolation  as  a  mathematical   theory  was  introduced  by  Broadbent  and
Hammersley \cite{broadbent:percolation},  as a stochastic way  of modeling the
flow of a fluid or gas through  a porous medium of small channels which may or
may not let gas or fluid pass.   It is one of the simplest models exhibiting a
phase transition,  and the occurrence of  a critical phenomenon  is central to
the appeal of percolation.  Having truly applied origins, percolation has been
used to model the fingering and spreading of oil in water, to estimate whether
one  can build  non-defective  integrated  circuits, to  model  the spread  of
infections and forest fires.  From a mathematical point of view percolation is
attractive   because   it  exhibits   relations   between  probabilistic   and
algebraic/topological properties of graphs.

\medskip

To make  the mathematical construction  of such a  system of channels,  take a
graph $\mathcal  G$ (which originally was  taken as $\ZZ^d$),  with vertex set
$\mathcal V$ and  edge set $\mathcal E$, and make  all the edges independently
\emph{open} (or  passable) with probability $p$ or  \emph{closed} (or blocked)
with probability $1-p$.  Write $P_p$ for the corresponding probability measure
on the set of configurations of open and closed edges --- that model is called
\emph{bond percolation}.   The collection  of open edges  thus forms  a random
subgraph of  $\mathcal G$, and the  original question stated  by Broadbent was
whether the  connected component of the  origin in that subgraph  is finite or
infinite.

A \emph{path} on $\mathcal G$ is a  sequence $v_1, v_2, \dots $ of vertices of
$\mathcal G$, such that for all  $i\ge 1$, $v_i$ and $v_{i+1}$ are adjacent on
$\mathcal G$.  A  path is called \emph{open} if all  the edges $\{v_i, v_{i+1}
\}$ between successive  vertices are open. The infiniteness  of the cluster of
the origin is  equivalent to the existence of an  unbounded open path starting
from the origin.

\medskip

There  is an  analogous model,  called \emph{site  percolation}, in  which all
edges are assumed  being passable, but the vertices  are independently open or
closed with  probability $p$ or  $1-p$, respectively. An  open path is  then a
path along which all vertices are open.  Site percolation is more general than
bond  percolation  in  the  sense  that  the existence  of  a  path  for  bond
percolation on a  graph $\mathcal G$ is equivalent to the  existence of a path
for site  percolation on  the covering graph  of $\mathcal G$.   However, site
percolation on a given graph may  not be equivalent to bond percolation on any
other graph.

\medskip

All graphs under consideration will be assumed to be connected, locally finite
and quasi-transitive.  If  $A, B \subset \mathcal V$,  then $A \leftrightarrow
B$ means that there exists an open path from some vertex of $A$ to some vertex
of $B$;  by a slight abuse of  notation, $u \leftrightarrow v$  will stand for
the existence  of a  path between  sites $u$ and  $v$, \emph{i.e.}   the event
$\{u\} \leftrightarrow  \{v\}$.  The \emph{open cluster} $C(v)$  of the vertex
$v$ is  the set of  all open vertices  which are connected  to $v$ by  an open
path:
$$C(v) = \{u \in {\mathcal V}: \; u \leftrightarrow v \}.$$
The  central  quantity of  the  percolation  theory  is the  \emph{percolation
  probability}:
$$\theta(p)  :=   P_p\{{\textrm  {\bf  0}}\leftrightarrow  \infty   \}  =  P_p
\{|C({\textrm {\bf 0}})|= \infty \}.$$

\bigskip

The most  important property of  the percolation model  is that it  exhibits a
\emph{phase transition},  \emph{i.e.} there exists a threshold  value $p_c \in
[0,1]$, such that the global behavior of the system is substantially different
in the  two regions $p <  p_c$ and $p >  p_c$.  To make  this precise, observe
that  $\theta$  is  a  non-decreasing   function.   This  can  be  seen  using
Hammersley's joint construction  of percolation systems for all  $p \in [0,1]$
on  $\mathcal G$:  Let $\{U(v),  v \in  \mathcal V  \}$ be  independent random
variables, uniform  in $[0,1]$.   Declare $v$ to  be $p$-open if  $U(v)\le p$,
otherwise it  is declared $p$-closed.  The configuration  of $p$-open vertices
has the distribution $P_p$ for each $p \in [0,1]$.  The collection of $p$-open
vertices is non-decreasing in $p$, and therefore $\theta(p)$ is non-decreasing
as well.  Clearly $\theta (0)= 0$ and $\theta (1) = 1$.

\begin{figure}[htbp]
  \centering
  \input{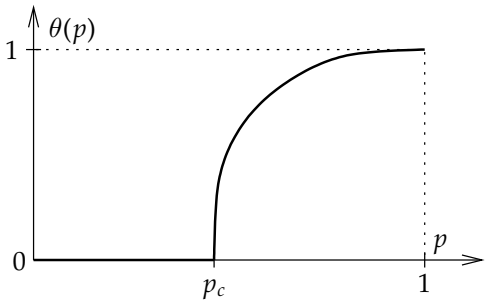}
  \caption{The  behavior  of  $\theta(p)$   around  the  critical  point  (for
    bond-percolation)}
  \label{fig:thetaofp}
\end{figure}

\bigskip

The \emph{critical probability} is defined as
$$p_c := p_c ({\mathcal G}) = \sup \{ p: \, \theta (p) = 0 \}.$$
By  definition, when $p<p_c$  the open  cluster of  the origin  is $P_p$-a.s.\
finite,  hence all  the clusters  are  also finite.   On the  other hand,  for
$p>p_c$ there is a strictly positive $P_p$-probability that the cluster of the
origin is infinite.  Thus, from Kolmogorov's zero-one law it follows that
$$P_p \{ |C(v)|  = \infty \; {\textrm {for some $v\in  {\mathcal V}$}}\!\} = 1
\;  {\textrm {for $p>p_c$}  }.$$ Therefore,  if the  intervals $[0,  p_c)$ and
$(p_c, 1]$ are both non-empty, there is a phase transition at $p_c$.

Using  a so-called  Peierls argument  it is  easy to  see  that $p_c({\mathcal
  G})>0$ for  any graph $\mathcal  G$ of bounded  degree.  On the  other hand,
Hammersley
proved that $p_c({\ZZ^d})  <1$ for bond percolation as soon  as $d\ge2$, and a
similar argument  works for  site percolation and  various periodic  graphs as
well.   But for  some graphs  $\mathcal G$  it  is not  so easy  to show  that
$p_c({\mathcal G})<1$.  One  says that the system is  in the \emph{subcritical
  (resp.\ supercritical) phase} if $p< p_c$ (resp.\ $p>p_c$).

It was one  of the most remarkable moments in the  history of percolation when
Kesten   proved   \cite{kesten:pc12}   that   the   critical   parameter   for
bond-percolation on $\ZZ^2$  is equal to $1/2$.  Nevertheless  the exact value
of $p_c({\mathcal  G})$ is  known only for  a handful  of graphs, all  of them
periodic and two-dimensional --- see below.

\section{Percolation in $\ZZ^d$}

The  graph on  which most  of the  theory was  originally built  is  the cubic
lattice $\ZZ^d$, and it was not  before the late 20th century that percolation
was seriously considered on other  kinds of graphs (such as \emph{e.g.} Cayley
graphs), on  which specific phenomena can  appear, such as  the coexistence of
multiple infinite clusters  for some values of the parameter  $p$. In all this
section,  the underlying  graph is  thus assumed  to be  $\ZZ^d$  for $d\ge2$,
although  most  of  the  results  still   hold  in  the  case  of  a  periodic
$d$-dimensional lattice.

\subsection{The sub-critical regime}


When $p<p_c$, all open clusters are  finite almost surely. One of the greatest
challenges  in percolation  theory  has been  to  prove that  $\chi (p):=  E_p
\{|C(v)|\}$  is finite  if  $p<p_c$  ($E_p$ stands  for  the expectation  with
respect to $P_p$). For that one can define another critical probability as the
threshold value  for the finiteness  of the expected  cluster size of  a fixed
vertex:
$$p_T({\mathcal G}) := \sup\{p:\, \chi (p) < \infty \}.$$
It was an  important step in the  development of the theory to  show that $p_T
({\mathcal  G})  = p_c  ({\mathcal  G})$.   The  fundamental estimate  in  the
subcritical regime,  which is a  much stronger statement than  $p_T ({\mathcal
  G}) = p_c ({\mathcal G})$, is the following:
\begin{thm}[Aizenman and Barsky, Menshikov]
  \label{t1} Assume  that $\mathcal G$ is  periodic. Then for $p  < p_c$ there
  exist constants $0 < C_1, C_2 < \infty$, such that
  $$P_p \{ |C(v)| \ge n\} \le C_1 e^{-C_2n}.$$
\end{thm}

The last statement can be sharpened to a ``local limit theorem'' with the help
of a  subadditivity argument :  For each $p<p_c$  there exists a  constant $0<
C_3(p) < \infty$, such that
$$\lim_{n \to \infty} - \frac{1}{n} \log P_p \{ |C(v)| = n\} = C_3 (p).$$

\subsection{The super-critical regime}

Once an infinite open cluster exists, it  is natural to ask how it looks like,
and  how many  infinite  open clusters  exist.   It was  shown  by Newman  and
Schulman that for periodic graphs, for  each $p$, exactly one of the following
three  situations  prevails: If  $N\in\ZZ_+\cup\{\infty\}$  is  the number  of
infinite open  clusters, then $P_p  (N=0) =  1$, or $P_p  (N=1) = 1$,  or $P_p
(N=\infty) =1$.


Aizenman,  Kesten and  Newman  showed that  the  third case  is impossible  on
$\ZZ^d$. By now  several proofs exist, perhaps the most  elegant proof of that
is due to  Burton and Keane, who prove that indeed  there cannot be infinitely
many infinite  open clusters on any  amenable graph.  However,  there are some
graphs,  such as  regular  trees,  on which  coexistence  of several  infinite
clusters is possible.

The geometry  of the infinite  open cluster can  be explored in some  depth by
studying the behavior of a random walk  on it. When $d=2$, the random walk is
recurrent, and when  $d \ge 3$ is a.s.\ transient. In  all dimensions $d\ge 2$
the walk behaves diffusively, and the Central Limit Theorem and the Invariance
principle were established in both the annealed and quenched cases.

\subsubsection*{Wulff droplets}

In  the  supercritical regime,  aside  from  the  infinite open  cluster,  the
configuration contains finite clusters  of arbitrary large sizes.  These large
finite  open clusters  can be  thought of  as droplets  swimming in  the areas
surrounded by an infinite open  cluster. The presence at a particular location
of a large finite cluster is  an event of low probability, namely, on $\ZZ^d$,
$d\ge 2$,  for $p>  p_c$, there exist  positive constants $0<  C_4(p), C_5(p)<
\infty$, such that
$$C_4(p) \leq  -\frac{1}{n^{(d-1)/d}} \log P_p  \{ |C(v)| = n\}  \leq C_5(p)$$
for all large $n$. This estimate is based on the fact that the occurrence of a
large finite cluster is due to a surface effect.  The typical structure of the
large finite cluster is described by the following theorem:

\begin{thm} 
  Let $d\ge2$, and $p>p_c$.  There exists a bounded, closed, convex subset $W$
  of $\mathbb  R^d$ containing the  origin, called the  \emph{normalized Wulff
    crystal}  of  the  Bernoulli  percolation  model,  such  that,  under  the
  conditional probability $P_p \{ \cdot  \mid n^d \le |C({\textrm {\bf 0}})| <
  \infty \}$, the random measure
  $$\frac{1}{n^d} \sum_{x \in C({\textrm {\bf 0}})} \delta_{x/n}$$  
  (where $\delta_{\mathbf{x}}$ denotes a Dirac mass at $\mathbf{x}$) converges
  weakly  in probability  towards the  random  measure $\theta(p)\un_W(x-M)dx$
  (where $M$ is  the rescaled center of mass  of the cluster $C(\mathbf{0})$).
  The deviation probabilities behave as $\exp\{-cn^{d-1} \}$ (\emph{i.e.} they
  exhibit large deviations of surface order). 
\end{thm}

This  result was  proved  in dimension  $2$  by Alexander,  Chayes and  Chayes
\cite{alexander:wulff},   and   in   dimensions   $3$   and   more   by   Cerf
\cite{cerf:wulff}.

\subsection{Percolation near the critical point}

\subsubsection{Percolation in slabs}

The  main  macroscopic observable  in  percolation  is  $\theta(p)$, which  is
positive    above    $p_c$,   $0$    below    $p_c$,    and   continuous    on
$[0,1]\setminus\{p_c\}$.   Continuity at  $p_c$  is an  open  question in  the
general case; it is  known to hold in two dimensions (cf.\  below) and in high
enough dimension  (at the  moment $d\geq19$ though  the value of  the critical
dimension is believed to be  $6$) using lace expansion methods. The conjecture
that $\theta(p_c)  = 0$ for  $3 \le d  \le 18$ remains  one of the  major open
problems.

Efforts to prove that led  to some interesting and important results.  Barsky,
Grimmett  and  Newman  solved  the   question  in  the  half-space  case,  and
simultaneously  showed that  the slab  percolation and  half-space percolation
thresholds coincide.  This was  complemented by Grimmett and Marstrand showing
that
$$p_c(slab)= p_c(\ZZ^d).$$

\subsubsection{Critical exponents}

In  the sub-critical regime,  exponential decay  of the  correlation indicates
that there  is a finite  \emph{correlation length} $\xi(p)$ associated  to the
system, and defined (up to constants) by the relation
$$P_p(0  \leftrightarrow   n\mathbf{x})  \approx   \exp  \left(  -   \frac  {n
    \varphi(x)} {\xi(p)} \right)$$
where  $\varphi$   is  bounded   on  the  unit   sphere  (this  is   known  as
\emph{Ornstein-Zernike decay}).  The phase transition can then also be defined
in terms  of the divergence  of the correlation  length, leading again  to the
same value for $p_c$;  the behavior at or near the critical  point then has no
finite   characteristic  length,   and   gives  rise   to  scaling   exponents
(conjecturally in most cases).

The most  usual critical exponents are  defined as follows,  if $\theta(p)$ is
the percolation probability,  $C$ the cluster of the  origin, and $\xi(p)$ the
correlation length:
\begin{align*}
  \frac{\partial^3}{\partial p^3} E_p[|C|^{-1}] &\approx |p-p_c|^{-1-\alpha}\\
  \theta(p) &\approx (p-p_c)_+^\beta\\
  \chi^f(p) := E_p[|C| \un_{|C|<\infty}] &\approx |p-p_c|^{-\gamma}\\
  P_{p_c}[|C|=n] &\approx n^{-1-1/\delta}\\
  P_{p_c}[x\in C] &\approx |x|^{2-d-\eta}\\
  \xi(p) &\approx |p-p_c|^\nu\\
  P_{p_c}[\mathrm{diam}(C)=n] &\approx n^{-1-1/\rho}\\
  \frac {E_p[|C|^{k+1} \un_{|C|<\infty}]} {E_p[|C|^k \un_{|C|<\infty}]} &\approx
  |p-p_c|^{-\Delta}
\end{align*}

These exponents are  all expected to be universal,  \emph{i.e.} to depend only
on the dimension  of the lattice, although this is not  well understood at the
mathematical  level;  the   following  \emph{scaling  relations}  between  the
exponents are believed to hold:
$$2-\alpha = \gamma+2\beta = \beta(\delta+1),\; \Delta = \delta \beta,\;
\gamma = \nu(2-\eta).$$
In addition,  in dimensions up  to $d_c=6$, two  additional \emph{hyperscaling
  relations} involving $d$ are strongly conjectured to hold:
$$d\rho = \delta+1,\quad d\nu=2-\alpha,$$
while above $d_c$  the exponents are believed to  take their mean-field value,
\emph{i.e.} the ones they have for percolation on a regular tree:
$$\alpha=-1,\; \beta=1,\; \gamma=1,\; \delta=2,$$
$$\eta=0,\;  \nu=\frac12,\;\rho=\frac12,\; \Delta=2.$$

\bigskip

Not much is  known rigorously on critical exponents in  the general case. Hara
and Slade  (\cite{hara-slade:meanfield}) proved that mean  field behavior does
happen above dimension $19$, and the proof can likely be extended to treat the
case  $d\geq7$.   In  the  two-dimensional  case on  the  other  hand,  Kesten
(\cite{kesten:scaling}) showed that, assuming  that the exponents $\delta$ and
$\rho$ exist, then so do $\beta$, $\gamma$, $\eta$ and $\nu$, and they satisfy
the scaling and hyperscaling relations where they appear.

\subsubsection{The incipient infinite cluster}

When studying long-range properties of a  critical model, it is useful to have
an  object which is  infinite at  criticality, and  such is  not the  case for
percolation  clusters. There  are two  ways to  condition the  cluster  of the
origin to be infinite  when $p=p_c$: The first one is to  condition it to have
diameter at  least $n$  (which happens with  positive probability) and  take a
limit in distribution  as $n$ goes to infinity; the second  one is to consider
the model for  parameter $p>p_c$, condition the cluster of  $0$ to be infinite
(which happens with positive probability)  and take a limit in distribution as
$p$ goes  to $p_c$. The limit  is the same in  both cases, it is  known as the
\emph{incipient infinite cluster}.

As  in  the  super-critical  regime,  the  structure of  the  cluster  can  be
investigated by studying the behavior of a random walk on it, as was suggested
by de  Gennes; Kesten proved  that in two  dimensions, the random walk  on the
incipient  infinite cluster  is  sub-diffusive, \emph{i.e.}   the mean  square
displacement  after   $n$  steps  behaves  as   $n^{1-\varepsilon}$  for  some
$\varepsilon>0$.

The construction of  the incipient infinite cluster was done  by Kesten in two
dimensions  \cite{kesten:iic},  and   a  similar  construction  was  performed
recently in high dimension by Van der Hofstad and Jarai (\cite{hofstad:iic}).

\section{Percolation in two dimensions}

As is  the case for several  other models of  statistical physics, percolation
exhibits  many  specific  properties  when  considered  on  a  two-dimensional
lattice: Duality arguments  allow for the computation of  $p_c$ in some cases,
and  for the  derivation  of \emph{a  priori}  bounds for  the probability  of
crossing  events at  or near  the  critical point,  leading to  the fact  that
$\theta(p_c)=0$.   On   another  front,   the   scaling   limit  of   critical
site-percolation on the two-dimensional triangular lattice can be described in
terms of SLE processes.

\begin{figure}[htbp]
  \centering
  \includegraphics[scale=0.33]{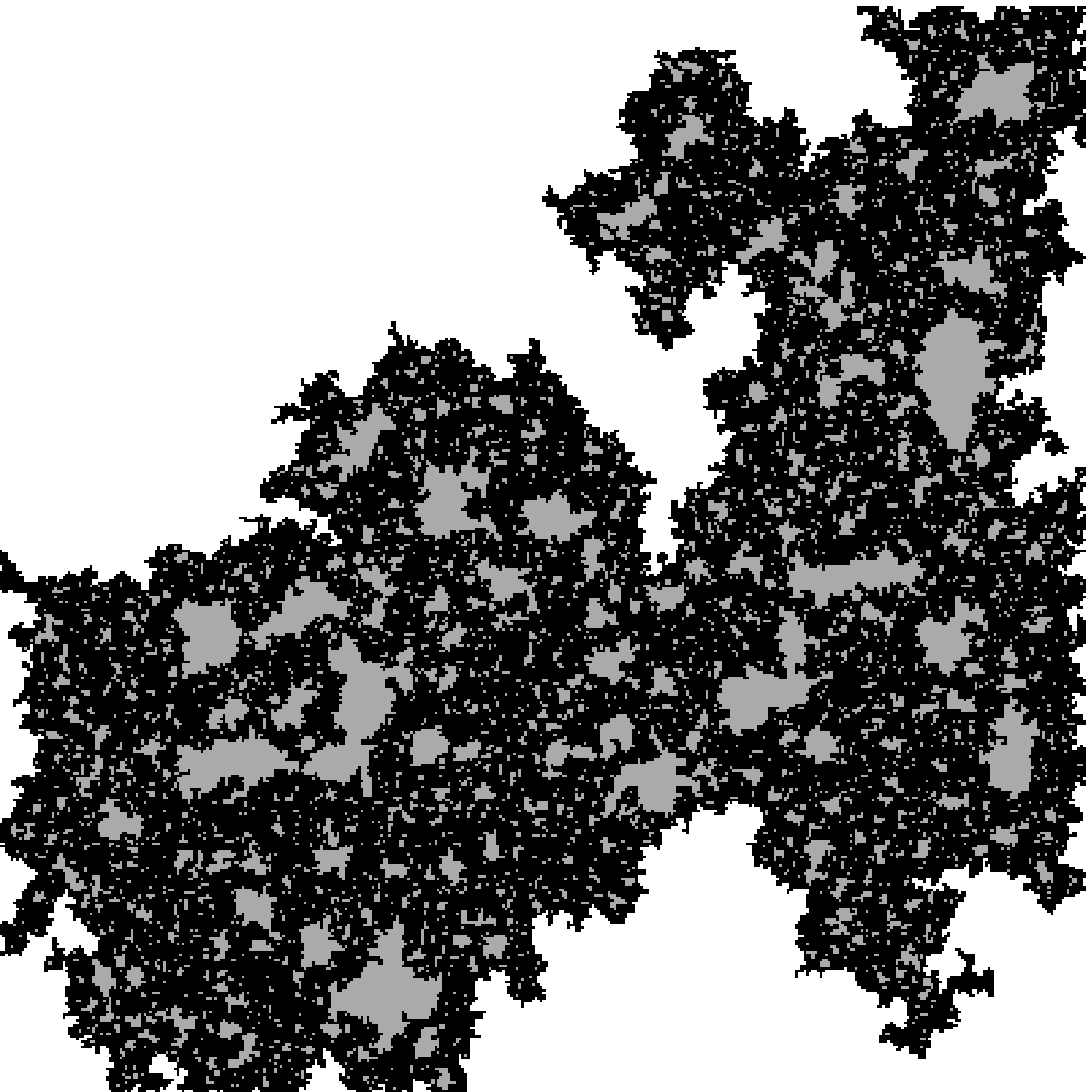}
  \hspace{1cm}
  \includegraphics[scale=0.33]{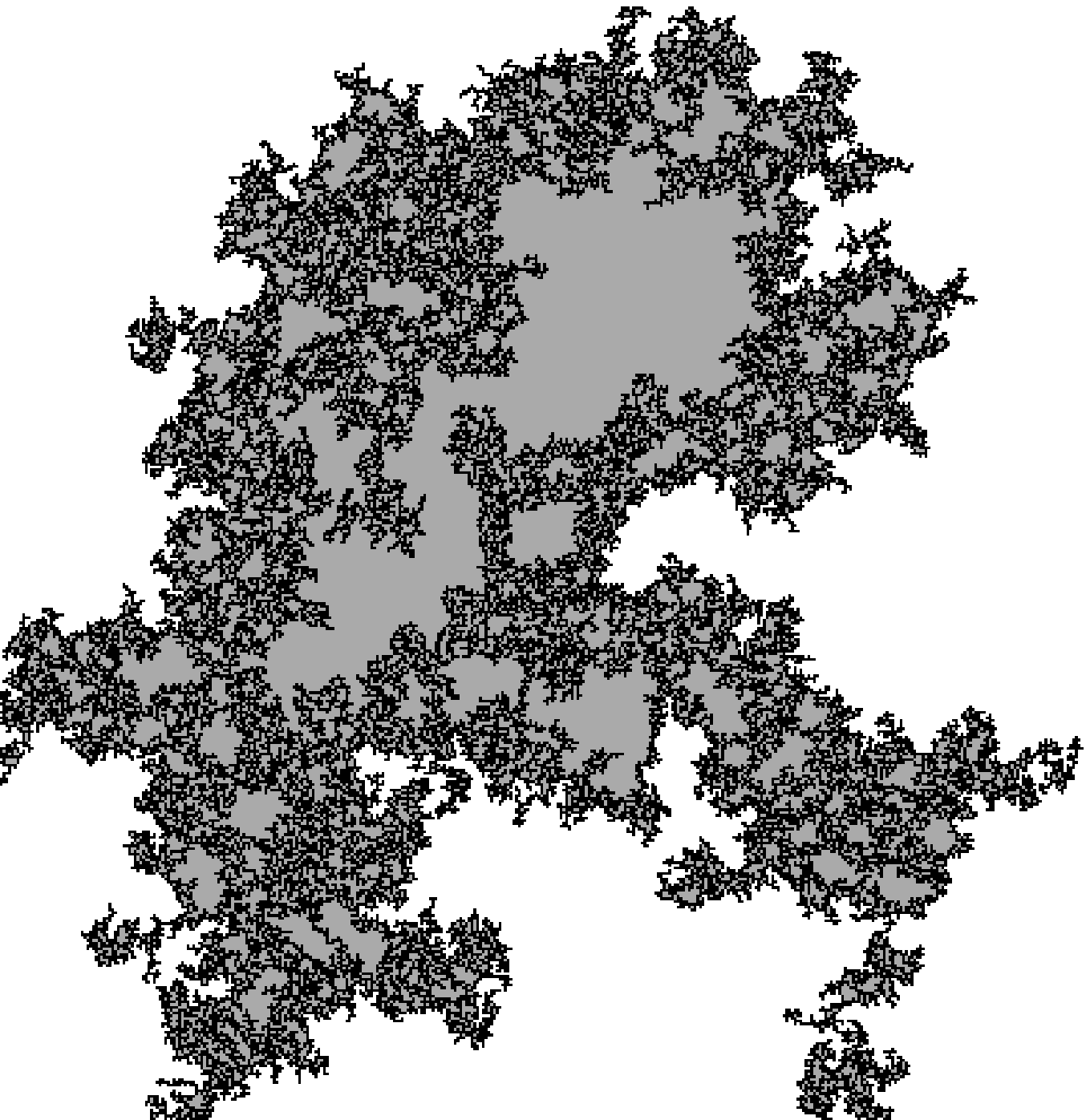}
  \caption{Two  large critical  percolation clusters  in a  box of  the square
  lattice (first: bond-percolation, second: site-percolation)}
  \label{fig:perc2d}
\end{figure}

\subsection{Duality, exact computations and RSW theory}

Given a planar lattice $\mathcal  L$, define two associated graphs as follows.
The \emph{dual  lattice} $\mathcal  L'$ has  one vertex for  each face  of the
original  lattice,  and an  edge  between  two vertices  if  and  only if  the
corresponding  faces of  $\mathcal L$  share an  edge.  The  \emph{star graph}
$\mathcal L^\ast$  is obtained by adding  to $\mathcal L$ an  edge between any
two vertices belonging to the same  face ($\mathcal L^{\ast}$ is not planar in
general;   $(\mathcal   L,\mathcal  L^{\ast})$   is   commonly   known  as   a
\emph{matching  pair}).  Then,  a result  of  Kesten is  that, under  suitable
technical conditions,
$$p_c^{\mathrm{bond}}(\mathcal   L)  +   p_c^{\mathrm{bond}}(\mathcal   L')  =
p_c^{\mathrm{site}}(\mathcal L) + p_c^{\mathrm{site}}(\mathcal L^\ast) = 1.$$

Two cases are  of particular importance: The lattice  $\ZZ^2$ is isomorphic to
its  dual; the  triangular lattice  $\mathcal  T$ is  its own  star graph.  It
follows that
$$p_c^{\mathrm{bond}}(\ZZ^2)  = p_c^{\mathrm{site}}(\mathcal  T)  = \frac12.$$
The   only   other   critical   parameters   that  are   known   exactly   are
$p_c^{\mathrm{bond}}(\mathcal     T)=2\sin(\pi/18)$     (and    hence     also
$p_c^{\mathrm{bond}}$ for $\mathcal  T'$, \emph{i.e.}  the hexagonal lattice),
and  $p_c^{\mathrm{bond}}$ for  the bow-tie  lattice which  is a  root  of the
equation  $p^5-6p^3+6p^2+p-1=0$.   The value  of  the  critical parameter  for
site-percolation on $\ZZ^2$ might on the other hand never be known, it is even
possible that it is ``just a number'' without any other signification.

\medskip

Still using duality, one can  prove that the probability, for bond-percolation
on  the square  lattice  with parameter  $p=1/2$,  that there  is a  connected
component crossing  an $(n+1)\times  n$ rectangle in  the longer  direction is
exactly equal to  $1/2$.  This and clever arguments  involving the symmetry of
the lattice lead to the following result, proved independently by Russo and by
Seymour and Welsh and known as the RSW theorem:

\begin{thm}[Russo~\cite{russo:rsw}; Seymour-Welsh~\cite{seymour:rsw}]
  For  every $a,b>0$  there exist  $\eta>0$ and  $n_0>0$ such  that  for every
  $n>n_0$, the  probability that  there is a  cluster crossing an  $\lfloor na
  \rfloor  \times \lfloor  nb \rfloor$  rectangle  in the  first direction  is
  greater than $\eta$.
\end{thm}

The  most direct consequence  of this  estimate is  that the  probability that
there is a cluster going around an annulus of a given modulus is bounded below
independently of the  size of the annulus; in  particular, almost surely there
is \emph{some}  annulus around  $0$ in  which this happens,  and that  is what
allows to prove that $\theta(p_c)=0$ for bond-percolation on $\ZZ^2$.

\subsection{The scaling limit}

RSW-type estimates  give positive evidence that  a scaling limit  of the model
should exist; it  is indeed essentially sufficient to  show convergence of the
crossing probabilities  to a  non-trivial limit as  $n$ goes to  infinity. The
limit, which  should depend only on  the ration $a/b$, was  predicted by Cardy
using conformal fields theory methods.  A most celebrated result of Smirnov is
the proof of Cardy's formula in the case of site-percolation on the triangular
lattice $\mathcal T$:

\begin{thm}[Smirnov \cite{smirnov:perco}]
  Let $\Omega$ be a simply connected domain of the plane with four points $a$,
  $b$, $c$, $d$ (in that order)  marked on its boundary. For every $\delta>0$,
  consider a  critical site-percolation model on the  intersection of $\Omega$
  with $\delta \mathcal T$ and let $f_\delta(ab,cd;\Omega)$ be the probability
  that it contains a cluster connecting the arcs $ab$ and $cd$.  Then:
  \begin{enumerate}
  \item   $f_\delta(ab,cd;\Omega)$   has   a  limit   $f_0(ab,cd;\Omega)$   as
    $\delta\to0$;
  \item The limit is conformally  invariant, in the following sense: If $\Phi$
    is   a    conformal   map   from    $\Omega$   to   some    other   domain
    $\Omega'=\Phi(\Omega)$, and maps $a$ to $a'$, $b$ to $b'$, $c$ to $c'$ and
    $d$ to $d'$, then $f_0(ab,cd;\Omega)=f_0(a'b',c'd';\Omega')$;
  \item In  the particular  case when $\Omega$  is an equilateral  triangle of
    side length $1$ and vertices $a$, $b$  and $c$, and if $d$ is on $(ca)$ at
    distance $x\in(0,1)$ from $c$, then $f_0(ab,cd;\Omega)=x$.
  \end{enumerate}
\end{thm}

Point 3.\ in  particular is essential since it allows  to compute the limiting
crossing probabilities  in any  conformal rectangle. In  the original  work of
Cardy, he made his prediction in the  case of a rectangle, for which the limit
involves hypergeometric  functions; the  remark that the  equilateral triangle
gives rise to nicer formulae is originally due to Carleson. 

\medskip

\begin{figure}[htbp]
  \centering
  \includegraphics[scale=0.45]{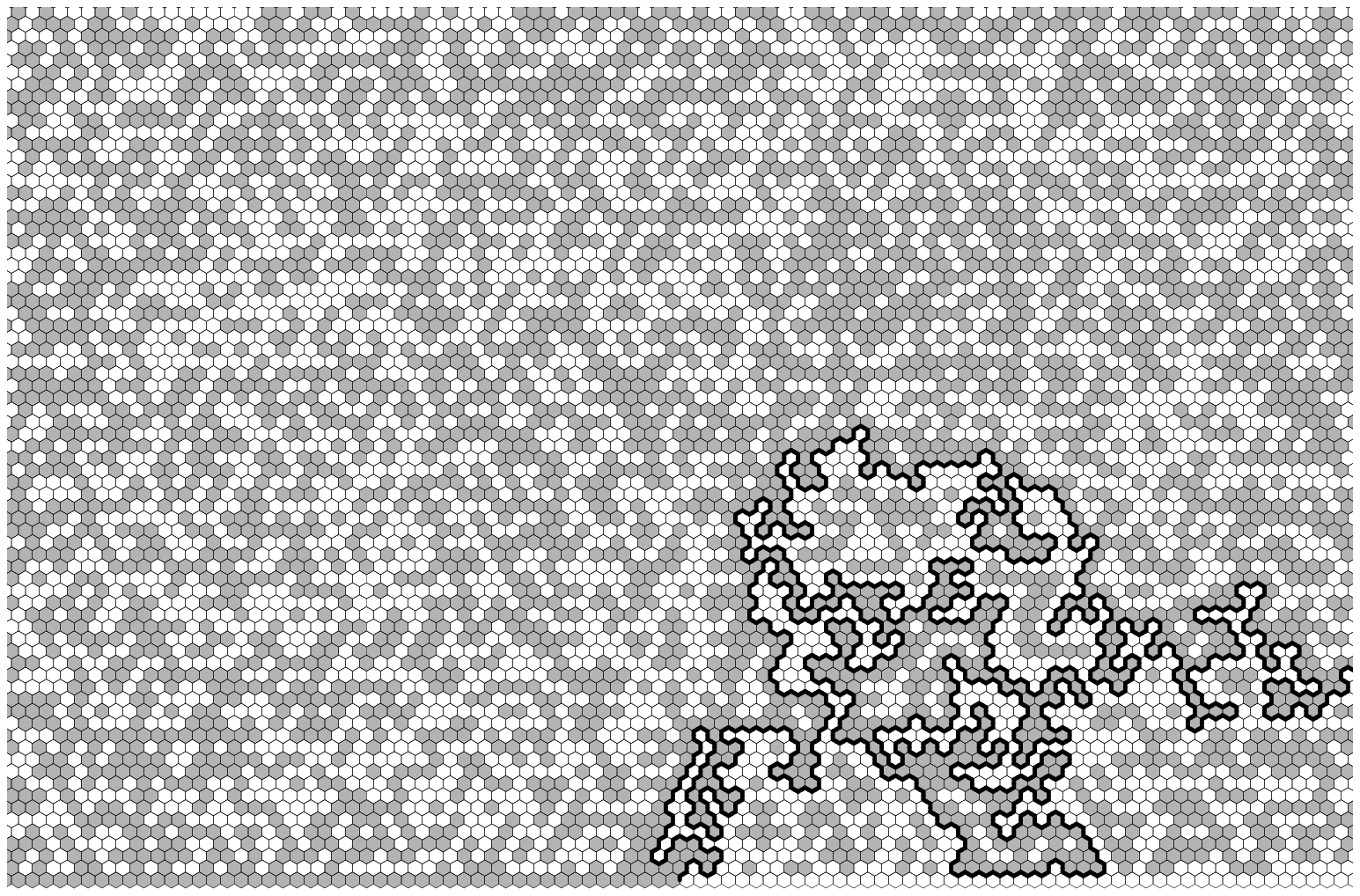}
  \caption{A percolation exploration path}
  \label{fig:expl}
\end{figure}

To precisely state the convergence of percolation to its scaling limit, define
the  random  curve  known  as  the \emph{percolation  exploration  path}  (see
fig.~\ref{fig:expl})  as   follows:  In  the  upper   half-plane,  consider  a
site-percolation model on  a portion of the triangular  lattice and impose the
boundary  conditions that on  the negative  real half-line  all the  sites are
open, while on the other half-line the sites are closed. The exploration curve
is then  the common boundary  of the open  cluster spanning from  the negative
half-line, and the closed cluster  spanning from the positive half-line; it is
an infinite, self-avoiding random curve in the upper half-plane.

\begin{figure}[htbp]
  \centering
  \includegraphics[scale=0.1]{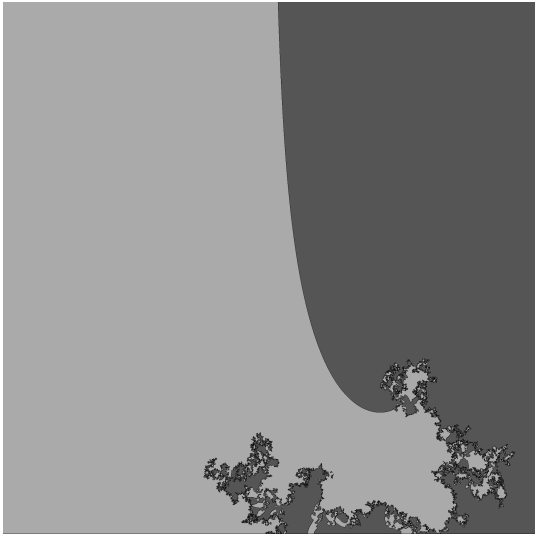}
  \caption{An $SLE$ process with parameter $\kappa=6$ (infinite time, with the
    driving process stopped at time $1$)}
  \label{fig:sle6}
\end{figure}

As the mesh  of the lattice goes to $0$, the  exploration curve then converges
in distribution to the trace of an $SLE$ process, as introduced by Schramm,
with parameter $\kappa=6$ ---  see fig.~\ref{fig:sle6}.  The limiting curve is
not simple  anymore (which  corresponds to the  existence of pivotal  sites on
large critical  percolation clusters), and  it has Hausdorff  dimension $7/4$.
For more details on $SLE$ processes,  see \emph{e.g.} the related entry in the
present volume.

\medskip

As an application of this convergence  result, one can prove that the critical
exponents described in the previous section do exist (still in the case of the
triangular  lattice), and  compute their  exact values,  except  for $\alpha$,
which is still listed here for completeness:
$$\left[\alpha=-\frac23,\right]\; \beta=\frac{5}{36},\;
\gamma=\frac{43}{18},\; \delta=\frac{91}{5},$$
$$\eta=\frac5{24},\;            \nu=\frac43,\;            \rho=\frac{48}{5},\;
\Delta=\frac{91}{36}.$$ 
These exponents  are expected to be  \emph{universal}, in the  sense that they
should be the same for percolation  on any two-dimensional lattice; but at the
time  of this  writing  this phenomenon  is  far from  being  understood on  a
mathematical level.

The rigorous  derivation of the critical  exponents for percolation  is due to
Smirnov and  Werner \cite{smirnov:exps}; the  dimension of the  limiting curve
was obtained by Beffara \cite{beffara:dimsle}.

\section{Other lattices and percolative systems}

Some  modifications or  generalizations of  standard Bernoulli  percolation on
$\ZZ^d$ exhibit an interesting behavior  and as such provide some insight into
the original  process as well; there  are too many  mathematical objects which
can be argued to be percolative in some sense to give a full account of all of
them, so the following list is somewhat arbitrary and by no means complete.

\subsection{Percolation on non-amenable graphs}

The  first modification  of  the  model one  can  think of  is  to modify  the
underlying graph and move away  from the cubic lattice; phase transition still
occurs,  and  the main  difference  is  the  possibility for  infinitely  many
infinite clusters  to coexist. On  a regular tree,  such is the  case whenever
$p\in(p_c,1)$,  the first  non-trivial example  was produced  by  Grimmett and
Newman as the  product of $\ZZ$ by a  tree: There, for some values  of $p$ the
infinite  cluster  is  unique,  while  for  others  there  is  coexistence  of
infinitely many  of them. The  corresponding definition, due to  Benjamini and
Schramm, is then the following: If $N$ is as above the number of infinite open
clusters,
$$p_u := \inf \left\{ p: P_p(N=1) = 1\right\} \ge p_c.$$ 
The main question is then to characterize graphs on which $0<p_c<p_u<1$.

A  wide  class  of interesting  graphs  is  that  of \emph{Cayley  graphs}  of
infinite,  finitely generated  groups.   There, by  a  simultaneous result  by
Häggström  and Peres  and by  Schonmann, for  every $p\in(p_c,p_u)$  there are
$P_p$-a.s.\ infinitely  many infinite  cluster, while for  every $p\in(p_u,1]$
there is only one --- note that this does not follow from the definition since
new infinite components could appear when $p$ is increased.  It is conjectured
that  $p_c<p_u$  for  any Cayley  graph  of  a  non-amenable group  (and  more
generally for any quasi-transitive  graph with positive Cheeger constant), and
a result by Pak and Smirnova
is that  every infinite, finitely  generated, non-amenable group has  a Cayley
graph on which $p_c>p_u$; this is then expected not to depend on the choice of
generators. In the general case, it was recently proved by Gaboriau
that if the  graph $\mathcal G$ is unimodular,  transitive, locally finite and
supports  non-constant  harmonic  Dirichlet functions  (\emph{i.e.}   harmonic
functions  whose   gradient  is  in  $\ell^2$),   then  indeed  $p_c({\mathcal
  G})<p_u({\mathcal G})$.

For reference and further reading on the topic, the reader is advised to refer
to  the review  paper by  Benjamini and  Schramm  \cite{benjamini:beyond}, the
lecture  notes of Peres\cite{peres:stflour},  and the  more recent  article of
Gaboriau \cite{gaboriau:percolation}.

\subsection{Gradient percolation}

\begin{figure}[htbp]
  \centering
  \includegraphics[scale=0.1]{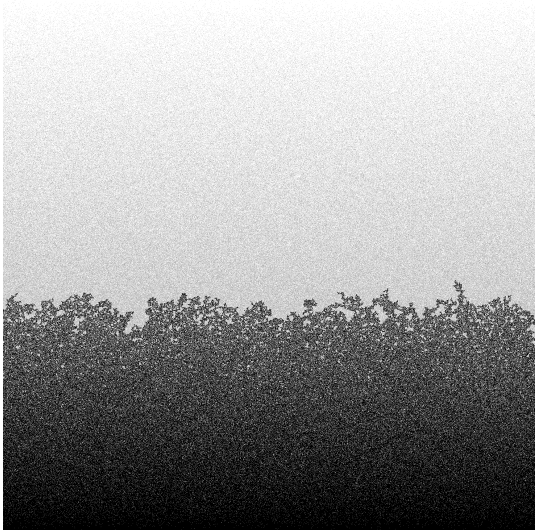}
  \caption{Gradient percolation in a square.  In black is the cluster spanning
  from the bottom side of the square.}
  \label{fig:grad}
\end{figure}

Another possible modification of the  original model is to allow the parameter
$p$ to  depend on the location; the  porous medium may for  instance have been
created by some kind of erosion, so  that there will be more open edges on one
side of a given  domain than on the other. If $p$  still varies smoothly, then
one expects some regions to look subcritical and others to look supercritical,
with  interesting  behavior  in  the   vicinity  of  the  critical  level  set
$\{p=p_c\}$.    This   particular  model   was   introduced   by  Sapoval   et
al.~\cite{sapoval:gradient} under the name of \emph{gradient percolation}; see
fig.~\ref{fig:grad}.

The control of  the model away from the critical zone  is essentially the same
as for  usual Bernoulli percolation, the  main question being  how to estimate
the width  of the  phase transition.   The main idea  is then  the same  as in
scaling theory: If the distance between a point $v$ and the critical level set
is less than  the correlation length for parameter $p(v)$, then  $v$ is in the
phase transition domain.  This of  course makes sense only asymptotically, say
in  a large  $n\times n$  square with  $p(x,y)=1-y/n$ as  is the  case  in the
figure: The transition then is expected  to have width of order $n^a$ for some
exponent $a>0$. 

\subsection{First-passage percolation}

First  passage  percolation (also  known  as  Eden  or Richardson  model)  was
introduced  by  Hammersley and  Welsh  in  1965  \cite{hammersley-welsh} as  a
time-dependent  model for the  passage of  fluid through  a porous  medium. To
define  the model, with  each edge  $e\in \mathcal  E(\ZZ^d)$ is  associated a
random variable  $T(e)$, which can be  interpreted as being  the time required
for fluid to flow along $e$. $T(e)$ are assumed to be independent non-negative
random variables having common distribution $F$.  For any path $\pi$ we define
the passage time $T(\pi)$ of $\pi$ as
$$T(\pi):= \sum_{e\in \pi} T(e).$$
The \emph{first passage  time} $a(x,y)$ between vertices $x$  and $y$ is given
by
$$a(x,y)=\inf \{T(\pi):{\textrm{$\pi$ a path from $x$ to $y$}}\};$$
and we can define
$$W(t):= \{ x \in \ZZ^d : a(0,x) \le t \},$$
the set  of vertices reached  by the  liquid by time  $t$.  It turns  out that
$W(t)$  grows approximately  linearly as  time passes,  and that  there exists
a non-random limit set $B$ such that either $B$ is compact and
$$(1-\epsilon)B \subseteq  \frac{1}{t}\widetilde W(t) \subseteq (1+\epsilon)B,
\;{\textrm{eventually a.s.}}$$
for all $\epsilon >0$, or $B= \Bbb R^d$, and
$$\{x\in \Bbb R^d  : \; |x|\le  K \}  \subseteq \frac{1}{t}\widetilde W(t), \;
{\textrm{eventually a.s.}}$$
for all $K>0$. Here $\widetilde W(t) = \{z + [-1/2, 1/2 ]^d : z \in W(t) \}$.

Studies  of first  passage percolation  brought many  fascinating discoveries,
including Kingman's celebrated sub-additive  ergodic theorem.  In recent years
interest  has been focused  on study  of fluctuations  of the  set $\widetilde
W(t)$  for large  $t$.   In spite  of  huge effort  and  some partial  results
achieved, it  still remains a  major task to establish  rigorously conjectures
predicted  by Kardar-Parisi-Zhang  theory  about shape  fluctuations in  first
passage percolation.

\subsection{Contact processes}

Introduced  by  Harris  and  conceived with  biological  interpretation,  the
\emph{contact process}  on $\ZZ^d$ is a continuous-time  process taking values
in the  space of subsets  of $\ZZ^d$. It  is informally described  as follows:
Particles are  distributed in $\ZZ^d$ in such  a way that each  site is either
empty or occupied  by one particle. The evolution  is Markovian: Each particle
disappears after an exponential time  of parameter $1$, independently from the
others; at any time, each particle  has a possibility to create a new particle
at any  of its  empty neighboring  sites, and does  so with  rate $\lambda>0$,
independently of everything else.


The question is  then whether, starting from a  finite population, the process
will die out  in finite time or whether it will  survive forever with positive
probability. The outcome will depend on the value of $\lambda$, and there is a
critical value $\lambda_c$, such that for $\lambda \le \lambda_c$ process dies
out, while  for $\lambda >  \lambda_c$ indeed there  is survival, and  in this
case  the shape of  the population  obey a  shape theorem  similar to  that of
first-passage percolation.

\medskip

The analogy with percolation  is strong, the corresponding percolative picture
being  the following:  In $\ZZ_+^{d+1}$,  each edge  is open  with probability
$p\in(0,1)$,   and  the  question   is  whether   there  exists   an  infinite
\emph{oriented} path  $\pi$ (\emph{i.e.}   a path along  which the sum  of the
coordinates is  increasing), composed of open  edges.  Once again,  there is a
critical parameter customarily  denoted by $\vec{p}_c$, at which  no such path
exists (compare  this to the open  question of the continuity  of the function
$\theta$ at $p_c$ in dimensions $3\le d\le 18$).
This variation of percolation lies  in a different universality class than the
usual Bernoulli model. 


\subsection{Invasion percolation}

Let $X(e): e \in {\mathcal E}$  be independent random variables indexed by the
edge set $\mathcal  E$ of $\ZZ^d$, $d\ge 2$,  each having uniform distribution
in  $[0,1]$. One  constructs a  sequence $C  = \{  C_i, i\ge  1\} $  of random
connected subgraphs of  the lattice in the following  iterative way: The graph
$C_0$ contains only the origin. Having defined $C_i$, one obtains $C_{i+1}$ by
adding to  $C_i$ an edge $e_{i+1}$  (with its outer  lying end-vertex), chosen
from the outer  edge boundary of $C_i$ so as  to minimize $X(e_{i+1})$.  Still
very little is known about the behavior of this process.


An interesting  observation, relating $\theta(p_c)$ of  usual percolation with
the invasion dynamics, comes from C.M. Newman:
$$\theta (p_c) = 0 \; \Leftrightarrow P \{x \in C \} \rightarrow 0 \; {\textrm
  {as }}\; |x| \rightarrow \infty.$$ 

\section*{Further reading}

For  a  much  more  in-depth   review  of  percolation  on  lattices  and  the
mathematical methods involved in its study,  and for the proofs of most of the
results we could  only point at, we  refer the reader to the  standard book of
Grimmett  \cite{grimmett:book}; another excellent  general reference,  and the
only place to  find some of the technical  graph-theoretical details involved,
is the  book of  Kesten \cite{kesten:book}.  More  information in the  case of
graphs  that are  not lattices  can be  found in  the lecture  notes  of Peres
\cite{peres:stflour}.

For curiosity, the reader can refer to the first mention of a problem close to
percolation,  in the  problem  section of  the  first volume  of the  American
Mathematical Monthly  \cite{wood:old}. References on more  specific topics are
given at the end of each section.


\section*{See also}

\textbf{Introductory  article: Statistical  mechanics; 2D  Ising  model; Wulff
  droplets; Stochastic Loewner evolutions.}

\bibliographystyle{siam}
\bibliography{perc}

\begin{thebibliography}{10}

\bibitem{alexander:wulff}
{\sc K.~Alexander, J.~T. Chayes, and L.~Chayes}, {\em The {W}ulff construction
  and asymptotics of the finite cluster distribution for two-dimensional
  {B}ernoulli percolation}, Comm. Math. Phys., 131 (1990), pp.~1--50.

\bibitem{beffara:dimsle}
{\sc V.~Beffara}, {\em Hausdorff dimensions for {$\rm SLE\sb 6$}}, Ann.
  Probab., 32 (2004), pp.~2606--2629.

\bibitem{benjamini:beyond}
{\sc I.~Benjamini and O.~Schramm}, {\em Percolation beyond {$\mathbb Z^d$},
  many questions and a few answers}, Electron. Comm. Probab., 1 (1996),
  pp.~no.\ 8, 71--82 (electronic).

\bibitem{broadbent:percolation}
{\sc S.~R. Broadbent and J.~M. Hammersley}, {\em Percolation processes, {I} and
  {II}}, Proc. Cambridge Philos. Soc., 53 (1957), pp.~629--645.

\bibitem{cerf:wulff}
{\sc R.~Cerf}, {\em Large deviations for three dimensional supercritical
  percolation}, vol.~267 of Ast\'erisque, SMF, 2000.

\bibitem{gaboriau:percolation}
{\sc D.~Gaboriau}, {\em Invariant percolation and harmonic {Dirichlet}
  functions}, Geometric And Functional Analysis,  (2005).
\newblock To appear.

\bibitem{grimmett:book}
{\sc G.~Grimmett}, {\em Percolation}, vol.~321 of Grundlehren der
  Mathematischen Wissenschaften, Springer-Verlag, Berlin, second~ed., 1999.

\bibitem{hammersley-welsh}
{\sc J.~M. Hammersley and D.~J.~A. Welsh}, {\em First-passage percolation,
  subadditive processes, stochastic networks, and generalized renewal theory},
  in Proc. Internat. Res. Semin., Statist. Lab., Univ. California, Berkeley,
  Calif., Springer-Verlag, New York, 1965, pp.~61--110.

\bibitem{hara-slade:meanfield}
{\sc T.~Hara and G.~Slade}, {\em Mean-field critical behaviour for percolation
  in high dimensions}, Communications in Mathematical Physics, 128 (1990),
  pp.~333--391.

\bibitem{kesten:pc12}
{\sc H.~Kesten}, {\em The critical probability of bond percolation on the
  square lattice equals {$1/2$}}, Comm. Math. Phys., 74 (1980), pp.~41--59.

\bibitem{kesten:book}
\leavevmode\vrule height 2pt depth -1.6pt width 23pt, {\em Percolation theory
  for mathematicians}, vol.~2 of Progress in Probability and Statistics,
  Birkh\"auser, Boston, Mass., 1982.

\bibitem{kesten:iic}
\leavevmode\vrule height 2pt depth -1.6pt width 23pt, {\em The incipient
  infinite cluster in two-dimensional percolation}, Probab. Theory Related
  Fields, 73 (1986), pp.~369--394.

\bibitem{kesten:scaling}
\leavevmode\vrule height 2pt depth -1.6pt width 23pt, {\em Scaling relations
  for {2D}-percolation}, Communications in Mathematical Physics, 109 (1987),
  pp.~109--156.

\bibitem{peres:stflour}
{\sc Y.~Peres}, {\em Probability on trees: an introductory climb}, in Lectures
  on probability theory and statistics (Saint-Flour, 1997), vol.~1717 of
  Lecture Notes in Math., Springer, Berlin, 1999, pp.~193--280.

\bibitem{russo:rsw}
{\sc L.~Russo}, {\em A note on percolation}, Z. Wahrscheinlichkeitstheorie und
  Verw. Gebiete, 43 (1978), pp.~39--48.

\bibitem{sapoval:gradient}
{\sc B.~Sapoval, M.~Rosso, and J.~Gouyet}, {\em The fractal nature of a
  diffusion front and the relation to percolation}, J. Phys. Lett., 46 (1985),
  pp.~L146--L156.

\bibitem{seymour:rsw}
{\sc P.~D. Seymour and D.~J.~A. Welsh}, {\em Percolation probabilities on the
  square lattice}, Ann. Discrete Math., 3 (1978), pp.~227--245.
\newblock Advances in graph theory (Cambridge Combinatorial Conf., Trinity
  College, Cambridge, 1977).

\bibitem{smirnov:perco}
{\sc S.~Smirnov}, {\em Critical percolation in the plane: Conformal invariance,
  {Cardy}'s formula, scaling limits}, C. R. Acad. Sci. Paris S\'er. I Math.,
  333 (2001), pp.~239--244.

\bibitem{smirnov:exps}
{\sc S.~Smirnov and W.~Werner}, {\em Critical exponents for two-dimensional
  percolation}, Mathematical Research Letters, 8 (2001), pp.~729--744.

\bibitem{hofstad:iic}
{\sc R.~van~der Hofstad and A.~A. J{\'a}rai}, {\em The incipient infinite
  cluster for high-dimensional unoriented percolation}, J. Statist. Phys., 114
  (2004), pp.~625--663.

\bibitem{wood:old}
{\sc D.~V. Wood}, {\em Average and probability, problem 5}, American
  Mathematical Monthly, 1 (1894), pp.~211--212.

\end{thebibliography}

\section*{Keywords}

Percolation,  random medium,  porous medium,  random graph,  phase transition,
critical  exponent,  shape   theorem,  correlation  length,  scaling,  scaling
relation, scaling limit, exploration path, {Cayley} graph, contact process.

\bigskip\bigskip

\begin{minipage}{5cm}
  \noindent\obeylines%
  Vincent BEFFARA
  UMPA -- ENS Lyon
  46 Allée d'Italie
  69364 Lyon Cedex 07
  FRANCE
\end{minipage}

\bigskip

\begin{minipage}{5cm}
  \noindent\obeylines%
  Vladas SIDORAVICIUS
  IMPA
  Estrada Dona Castorina 110
  Rio de Janeiro 22460-320
  BRASIL
\end{minipage}

\end{document}